\documentstyle[12pt]{article}
\newcommand{\vs}{\vspace{.2in}}

\newcommand{\noin}{\noindent}

\def \L{{ \it L}}
\def\qed{\hskip .6em \raise1.8pt\hbox{\vrule height4pt width6pt depth2pt}}

\def\R{\rm I\kern-2ptR}
\def\N{\rm I\kern-2ptN}
\def \supp{{\rm supp}}

\begin{document}

\begin{center}
{\bf A  rearrangement invariant space isometric to $L_p$ coincides with $L_p$}
\end{center}
\vs

\centerline {\bf Yuri Abramovich and Mikhail Zaidenberg}
\vs

\centerline{ Department of Mathematical Sciences, IUPUI,
Indianapolis, IN 46202, USA}

\centerline{ Institut Fourier des Math\'ematiques, Universit\'e Grenoble I,
BP 74, France}
\vs

\hfill {\it Dedicated to the memory of Ioseph Shneiberg}
\vs

The following theorem  is the main result of this note.\vs

\noin {\bf Theorem 1.}
{\it Let $(E, \|\cdot\|_E) $ be a rearrangement invariant Banach
function space on the interval $[0, 1]$. If $E$ is isometric to
$\L_p [0, 1]$ for some $1\le p<\infty$, then $E$ coincides with
$\L_p [0, 1]$ and furthermore $\|\cdot\|_E = \lambda\|\cdot\|_{\L_p}$,
where $\lambda = \|{\bf 1}\|_E$.}

\vs

We precede  the proof with some necessary definitions and notation. Two
Banach lattices $X$ and $Y$  are said to be {\bf order isometric} if there
exists an isometry $U$ of $X$ onto $Y$ which preserves the order, that is,
$U$ is an isometric surjective operator and $U(x) \ge 0$  if and only if $x
\ge 0$.\footnote{It is worth noticing that, as shown in [A1], a mere
positivity of an isometry $U$ implies that $U$ preserves the order.}

Let $L_0=L_0{[0,1]}$ be the vector lattice of all (equivalence classes of)
measurable real valued  functions on $[0,1]$ and let  $\mu$ denote
Lebesgue measure.

A Banach space $E$  is called  {\bf  rearrangement invariant} (r.i.)
if  the following three conditions hold:

\begin{itemize}
\item    $E$ is an ideal in $L_0$, i.e., if $x\in E,\ y\in L_0$
             and $|y| \le |x|$, then $y\in E$.

\item  If $x,y\in E$ and $|y| \le |x|$, then $\|y\| \le \|x\|$.

\item  If $x \in E$, $y \in L_0$ and the
            functions $|x|,\,|y|$ are equimeasurable\footnote {That is, $\mu
            \{t:|x(t)| < s\} = \mu  \{t: |y(t)| < s\}$
            for each $s \in \R $}, then $y \in E$ and $\|x\| = \|y\|$.
\end{itemize}

\noindent As usual, the symbol $E^+$ denotes the cone of all
nonnegative elements in $E$, and $E^{++}$ is a subset of all week
units\footnote {An element $0<x\in E$ is said to be a weak unit if
$x\land y >0$ for each $0<y\in E$.} of $E^+$. Recall that an ideal $B$ in $E$
is called a {\it band} if whenever $x=\sup x_\alpha $ for $x_\alpha \in B$ and
$x\in E$ we have $x\in B$. For each band $B$ in $E$ we denote by \supp~$B$
its support set, i.e., a unique (modulo subsets of measure zero) minimal
measurable subset $e \subseteq [0,1]$ such that $x\chi_e = x$ for each $x\in
B$, where $\chi_e$ is the characteristic function of the set $e$.

If $X$ and $Y$ are Banach lattices then, as usual, $X \oplus_p Y$ denotes
the Banach  lattice of all the pairs $(x, y) \in X \times Y$ with the norm
$\|(x, y)\| =(\|x\|^p + \|y\|^p )^{1\over p}$.

A key element  in our proof  is provided  by a theorem due to E.~Lacey and
P.~Wojtaszczyk [LW]. %We need  some notation from [LW].
Let $\l_p (2)$ denote
the standard two-dimensional $\l_p$-space, ordered by
the usual positive cone $\l_p^+ (2) = \{\alpha {\bf e}_1 + \beta {\bf e}_2 \,
|\,\alpha \ge 0,\,\beta\ge 0\}$, where ${\bf e}_1 = (1,0)\,,\,{\bf
e}_2 = (0,1)$.
Following [LW] we denote by $E_p (2)$ the same Banach space
$\l_p (2)$ but ordered by the cone
$$
E_p^+ (2) =  \{\alpha {\bar e}_1 + \beta {\bar e}_2\colon \
|\beta| \le \alpha \}.
$$
Both  $\l_p (2)$ and $E_p(2)$  are Banach lattices,  but they are not order
isometric. Let $\L_p (E_p (2))$ be the Banach space of all
(equivalence classes of) measurable vector-valued functions $f(t)$ on the
interval $[0,1]$ with values in $E_p(2)$, and with the norm  $\|f\| =
(\int\limits_0^1 \|f(t)\|^p\,dt)^{1\over p}$. It is plain to see
that under the natural ordering (i.e., $f\ge 0$ if and only if
$f(t) \in E_p^+(2)$) the space $\L_p (E_p (2))$ is a Banach lattice.

\vs
\noindent {\bf Theorem 2} (Lacey and Wojtaszczyk).
{\it  If a Banach lattice $E$ is  isometric to
$L_p[0,1] \ \ (1\le p\neq 2<\infty)$, then it is order
isometric to one of the next three Banach lattices:

{\rm (1)}  $\L_p [0,1]$ with the standard order;

{\rm (2)} $\L_p (E_p (2))$;

{\rm (3)}  $\L_p [0,1] \oplus_p \L_p (E_p (2))$.}

\vs

\noindent We point out that no two of  these  three Banach lattices are
order isometric. Finally, for each $t \in [0,1]$  we denote by  $\L_p (E_p(2);
[0,t])$ the  band in $\L_p (E_p (2))$ whose support set is  $[0,t]$.

The
unexplained terminology and notation regarding rearrangement invariant
spaces can be found in [KPS] and [LT], regarding Banach lattices in
[LT] and [V].

\vs

\noindent {\bf Proof of Theorem 1.} The case $p=2$  was proved in [S],  so in
what follows  we assume that $p \neq 2$.

Let $T$ be an isometry of the space $E$ onto
$\L_p [0,1]$. Therefore $T$ induces a new Banach lattice structure on the
Banach space $\L_p [0,1]$. By Theorem 2, this new Banach lattice is order
isometric to one of the three Banach lattices (1--3) indicated above.
We begin by showing   that Cases (2) and (3) are
impossible under the conditions of our theorem.

Suppose that Case (3) takes place. Therefore there exists an order isometry
$U: \L_p [0,1] \oplus_p \L_p (E_p (2)) \to E$.
Hence $U(\L_p [0,1])$
and $U( \L_p (E_p (2)))$ are complemented bands in $E$. We claim
the existence of the real numbers $t, \tau \in (0,1]$ such that
$$
\mu
(\supp \, U(\L_p [0,t])) = \mu (\supp \, U(\L_p (E_p (2); [0,
\tau]))) \eqno{(*)}
$$
Assume, for instance, that
$$
\mu (\supp \, U(\L_p [0,1])) = \alpha \ge \beta =  \mu
(\supp \, U(\L_p (E_p (2))).
$$
Put $A_t = \supp \, U(\L_p [0,t]))$ for $t
\in (0,1]$. A straightforward verification shows that the function
$t\mapsto \mu (A_t)$ is continuous.
Since
$\mu (A_1) =\alpha \ge \beta$ and $A_t \downarrow \emptyset$ when $t \to 0$,
we can conclude that $\mu (A_{t_0}) = \beta$ for some $t_0 \in (0,1]$.

The case when $\alpha < \beta$ can be treated similarly.
So in our  r.i. space   $E$, we have found two non-zero
bands $U(\L_p [0,t])$ and $U(\L_p (E_p (2); [0, \tau]))$, with supports of
equal measures. However, {\it in an arbitrary r.i. space any two
bands whose supports  have equal measures are clearly order isometric}.
A contradiction since the Banach lattices (1) and (2) are not order isometric.

Similar arguments can be applied to exclude Case (2). Indeed
for any  $0 < t \le 1/2$
the Banach lattice  $\L_p (E_p (2))$ contains two disjoint bands,
one of which is  $\L_p (E_p (2); [0, t])$
and another  is (order isometric to) $\L_p [0,t]$. Again this is
impossible as $E$ is a r.i. space.

Thus, we have established that the only case which can occur is the first one.
This  means that the r.i. spaces  $E$ and $\L_p [0,1]$ are isometric if
and only if they  are order isometric. However, a well known theorem due to
L.~Potepun [P1] (see also [A2] for a very simple proof of Potepun's theorem)
asserts that if two r.i. spaces $X$ and $Y$ (on the same measure space) are
order isomorphic, then they coincide. Therefore $E$ and $\L_p [0,1]$ coincide,
and consequently their norms are equivalent. To show that in actuality
$\|\cdot\|_E = \lambda \|\cdot\|_{\L_p}$ one can either make use of the
description (due to S. Banach [B]) of the isometry group of the space
$L_p[0,1]$, or (to bypass any technical details) simply to apply the
following theorem  due to D. A. Vladimirov [VS]:

\vs

\noindent {\bf Theorem 3.}
{\it Let $X$ be a r.i. KB-space. If the group of
the order isometries of $X$ acts transitively  on both
$X^{++}$ and $X^+\setminus X^{++}$, then $X$ is $\L_q [0,1]$ for some $q\ge 1$
and  $\|\cdot\|_X = \lambda \|\cdot\|_{\L_q}$ for some $\lambda > 0$. }
\vs

Since the space $ \L_p [0, 1]$ satisfies the assumptions of Theorem 3,
the space $E$, which is order isometric to  $\L_p [0, 1]$, satisfies them too.
So, by  Theorem 3,  $E = \L_q [0,1]$ for some $q\ge 1$ and
$\|\cdot\|_E = \lambda \|\cdot\|_{\L_q}$ for some $\lambda > 0$. Clearly,
$p=q$ and $\lambda = \|{\bf 1}\|_E$. This completes the proof
of Theorem~1. \qed
\vs

\noindent {\bf Remark.} If $p=\infty$, i.e., if  the r.i.\, space $E$ is
isometric to $L_\infty[0,1]$, then  the conclusion of Theorem~1 still
holds.

We will sketch the proof of the equality  $E= L_\infty[0,1]$.
It is well known (and  obvious) that an arbitrary r.i. space on $[0,1]$
contains the constant-one-function $\bf 1$, and thus $ L_\infty[0,1] \subseteq
E$.

Notice (see [AW, p. 324]) that  Banach lattice $E$ is
isomorphically  an AM-space, and consequently there exists a constant $M>0$ such
that  $\|e_1\vee e_2\vee \ldots\vee  e_n\|_E \le M$
for arbitrary positive elements $e_1,\ldots,e_n$ in the unit ball of $E$.
We are ready to verify now that the converse inclusion
$E\subseteq L_\infty[0,1]$
also holds, i.e., each function $x$ from $E_+$ is  (essentially) bounded.
Indeed, assume  contrary to what we claim,  that $x$ is an unbounded function.
Then for each $K>0$ there exists a set $A$ (depending on $x$ and $K$) of
positive measure such that $x(t) > K$ for almost all $t\in A$. Now using the
symmetry of the space $E$ we can find a finite number of functions  in $E$
each of which is equimeasurable with $x\chi_{A}$ and whose supremum $y$ is
greater than  function $K\bf 1$. This would imply that
$M\|x\|_E\ge \|y\|_E \ge K\|{\bf 1}\|_E$.
A contradiction, as $K$ is arbitrary. Thus, we have
established that  $E$  (as a set) coincides with $ L_\infty [0,1]$.
\qed

\vs

Some comments are in order in connection with  this article. Theorem 1
was proved by the authors many years ago. It was announced in [AZ],
but the proof was never published due to several reasons.
Our decision to publish it now has been inspired
by  a resent revival of interest in rearrangement invariant spaces.
This revival is due, in particular,
to a remarkable result by N. Kalton and  B.~Randrianantoanina [KR1,2] who
solved a longstanding problem   by showing that the description of
surjective isometries obtained for  the complex rearrangement invariant
spaces by the second author [Z1,2] (see also [Z3]) remains valid for the real
spaces as well. The proof of Theorem 1 presented above is basically our
original proof, it is very simple and is independent of a powerful
machinery developed later on in [JMST],
%It can be simplified substantially
%if we use either Kalton--Randrianantoanina's theorem (to show that the
%existence of an isometry implies the existence of an order isometry),
%or other relevant results from [JMST],
[K] and [LT].

We conclude with one remark and an open question. B.~Randrianantoanina
has obtained recently an analogue of Theorem 1 for the  Orlicz and
Lorentz spaces. An isomorphic version of the problem at hand
reads as follows:
\vs

\noindent {\bf  Problem.} {\em Let $E_1$, $E_2$ be two isomorphic r.i. spaces
on a measure space  $(\Omega, \Sigma, \lambda)$. When does it imply that the
spaces $E_1$ and $E_2$ coincide?}

In this case, of course, the norms will only be  equivalent.
When $E_2=L_2(\Omega, \Sigma, \lambda)$ this question was answered
in the  affirmative by  L.~Potepun [P2]. We refer to [JMST] and [K]
for several important results  when the answer to this problem is ``yes"
and for examples when the answer is ``no."

In the course of our work we had stimulating discussions with
V.~Ovchinnikov and L.~Yanovskii; it is a great pleasure  to thank them.
We also thank B.~Randrianantoanina for her interest in this work.
Last but not least we thank the referee who informed us about
a paper by  M. Hudzik, W. Kurc and M.Wis\l a [HKW], in which the Orlicz spaces
isometric to $L_\infty$ are characterized.

\vs

\centerline {\bf References}
\vs

\noin [A1]  Y. A. Abramovich, Some results on the isometries in normed
lattices,  {\em Optimization}, Novosibirsk $\,$ {\bf 43}(1988), 74--80.\\

\noin [A2]  Y. A. Abramovich, Operators preserving disjointness on
rearrangement invariant spaces, {\em Pacific. J. Math. \/ {\bf 148}}(1991),
201--206.\\

\noin [AW] Y. Abramovich and P. Wojtaszczyk,  The uniqueness of order
in the spaces $ L_p(0,1) $ and $ \ell_p  \/ (1\leq p \leq \infty ) $,
{\em Mat. Zametki\/ \bf18}$\,$(1975), 313--325.\\

\noin [AZ] Y. Abramovich, M. Zaidenberg, {\em Rearrangement invariant
spaces},   Notices AMS, {\bf 26} (1979), no. 2., A-207.\\

\noin [B] S. Banach, {\em Th\'eorie des op\'erations lin\'eaires},
Hafner Publishing Co, N.Y. 1932.  \\

\noin [HKW] M. Hudzik, W. Kurc and M.Wis\l a,  Strongly extreme points in Orlicz
function spaces, {\em J. Math. Anal. Appl.}, forthcoming.\\

\noin [JMST] W. B. Johnson, B. Maurey, G. Schechtman and L. Tzafriri, {\it
Symmetric Structures in Banach Spaces}, Mem. Amer. Math.Soc. No. 217, 1979.\\

\noin [K] N. J. Kalton, {\it Lattice Structures on Banach Spaces},
Mem. Amer. Math.Soc. No. 493, 1993.\\

\noin [KR1] N. J. Kalton, B. Randrianantoanina,  Isometries on
rearrangement-invariant spaces, {\em C. R. Acad. Sci. Paris} {\bf 316}
(1993), 351--355.\\

\noin [KR2] N. J. Kalton, B. Randrianantoanina,  Surjective
isometries on  rearrangement-invariant spaces, {\em Quart. J. Math.
Oxford}, forthcoming. \\

\noin [KPS] S. G. Krein, Y. I. Petunin, E. M. Semenov, {\em Interpolation
of  linear operators}, Providence, RI, 1982.\\

\noin [LW] E. Lacey and P. Wojtaszczyk,  Banach lattice structure on separable
$L_p$-spaces, {\em Proc. Amer. Math. Soc.\/ \bf 54}(1976), 83--89.\\

\noin [LT] J. Lindenstrauss and L. Tzafriri, {\em Classical Banach
Spaces {\rm II}, Function Spaces},  Springer Verlag,  Berlin--New York 1979.\\

\vfill
\eject

\noin [P1] L. I. Potepun, On isomorphisms of rearrangement invariant spaces,
          {\em Sibirsk. Math. J. \/ {\bf 12}}(1971), 623--629. \\

\noin [P2] L. I. Potepun, Hilbertian symmetric spaces, {\em Izv. VUZ-ov.
           Matematika}, {\bf 1 (140)} (1974), 90--95. \\

\noin [S] E. M. Semenov, {\em Interpolation of linear operators in
symmetric  spaces}, Thesis, Voronezh, 1968. \\

\noin [VS] D. A. Vladimirov and Yu. V. Shergin,  Characterization of
$L^p$-spaces in the class of normed vector lattices by means of their
isometric order automorphism group, {\em Vestnik LGU}, no. 1(1991),
17--20.\\

\noin [V] B. Z. Vulikh, {\em Introduction to the theory of partially
ordered  spaces}, Wolter--Noordhoff Sci. Publ., Ltd., Gr\"oningen, 1967.\\

\noin [Z1] M. G. Zaidenberg,  On the isometric classification of symmetric
              spaces, {\em Soviet Math. Dokl.} {\bf 18} (1977), 636--640.\\

\noin [Z2]  M. G. Zaidenberg,  Special representations of isometries
           of  functional spaces, {\em Investigations on the theory of
           functions  of several real variables}, Yaroslavl' 1980, 84--91.\\

\noin [Z3]  M. G. Zaidenberg, Groups of isometries of Orlich spaces,
        {\em Soviet   Math.  Dokl.} { \bf  17} (1976), 432--436.\\
%[2ex]

\end{document}